\title{Teaching Differentiation: A Rare Case for the Problem of the Slope of the Tangent Line}

\author{
Roman Kvasov\\
Department of Mathematics\\
University of Puerto Rico at Aguadilla \\
Aguadilla, Puerto Rico 00604, USA
}
\date{\today}

\documentclass[12pt]{article}
\usepackage{geometry}
\usepackage{setspace}
\usepackage{graphicx}
\usepackage{amssymb}
\usepackage{amsfonts}
\usepackage{amsmath}
\usepackage{array}
\usepackage{booktabs,tabularx}
\usepackage{float}
\usepackage{multicol}
\usepackage[utf8]{inputenc}
\usepackage[T1]{fontenc}

\begin{document}

\maketitle

\begin{abstract}
In this article we discuss an important students' misconception about derivatives, that the expression of the derivative of the function contains the information as to whether the function is differentiable or not where the expression is undefined. As a working example we consider a typical Calculus problem of finding the horizontal tangent lines of a function. Following the standard procedure, we derive the expression for the derivative using Product Rule. The search for the values of the independent variable, that make the derivative equal zero, leads to missing the unique solution of the problem. We show that in this case, even though the expression of the derivative is undefined, the function indeed possesses the derivative at the point. We also provide the methodological treatment of such functions, which can be effectively used in the classroom. 
\end{abstract}

\section{Introduction}

Derivate is one of the most important topics not only in mathematics, but also in physics, chemistry, economics and engineering. Every standard Calculus course provides a variety of exercises for the students to learn how to apply the concept of derivative. The types of problems range from finding an equation of the tangent line to the application of differentials and advanced curve sketching. Usually, these exercises heavily rely on such differentiation techniques as Product, Quotient and Chain Rules, Implicit and Logarithmic Differentiation \cite{Stewart2012}. The definition of the derivative is hardly ever applied after the first few classes and its use is not much motivated.

Like many other topics in undergraduate mathematics, derivative gave rise to many misconceptions \cite{Muzangwa2012}, \cite{Gur2007}, \cite{Li2006}. Just when the students seem to learn how to use the differentiation rules for most essential functions, the application of the derivative brings new issues. A common students' error of determining the domain of the derivative from its formula is discussed in \cite{Rivera2013} and some interesting examples of the derivatives, defined at the points where the functions themselves are undefined, are provided. However, the hunt for misconceptions takes another twist for the derivatives undefined at the points where the functions are in fact defined.

The expression of the derivative of the function obtained using differentiation techniques does not necessarily contain the information about the existence or the value of the derivative at the points, where the expression for the derivative is undefined. In this article we discuss a type of continuous functions that have the expression for the derivative undefined at a certain point, while the derivative itself at that point exists. We show, how relying on the formula for the derivative for finding the horizontal tangent line of a function, leads to a false conclusion and consequently to missing a solution. We also provide a simple methodological treatment of similar functions suitable for the classroom.

\section{Calculating the Derivative}

In order to illustrate how deceitful the expression of the derivative can be to a students' eye, let us consider the following problem.

\vspace{12pt}

\fbox{\begin{minipage}{5.25in}

\begin{center}

\begin{minipage}{5.0in}

\vspace{10pt}

\emph{Problem}

\vspace{10pt}

Differentiate the function $f\left(x\right)=\sqrt[3]{x}\sin{\left(x^2\right)}$. For which values of $x$ from the interval $\left[-1,1\right]$ does the graph of $f\left(x\right)$ have a horizontal tangent?

\vspace{10pt}

\end{minipage}

\end{center}

\end{minipage}}

\vspace{12pt}

Problems with similar formulations can be found in many Calculus books \cite{Stewart2012}, \cite{Larson2010}, \cite{Thomas2009}. Following the common procedure, let us find the expression for the derivative of the function $f\left(x\right)$ applying the Product Rule:
\begin{eqnarray}
f'\left(x\right) &=& \left(\sqrt[3]{x}\right)'\sin{\left(x^2\right)}+\left(\sin{\left(x^2\right)}\right)'\sqrt[3]{x} \notag \\ &=& \frac{1}{3\sqrt[3]{x^2}}\sin{\left(x^2\right)}+2x\cos{\left(x^2\right)}\sqrt[3]{x} \notag \\ &=& \frac{6x^2\cos{x^2}+\sin{x^2}}{3\sqrt[3]{x^2}} \label{DerivativeExpression}
\end{eqnarray}

Similar to \cite{Stewart2012}, we find the values of $x$ where the derivative $f'\left(x\right)$ is equal to zero:
\begin{equation}
6x^2\cos{x^2}+\sin{x^2} = 0 
\label{DerivativeEqualZero}
\end{equation}

Since the expression for the derivative (\ref{DerivativeExpression}) is not defined at $x=0$, it is not hard to see that for all values of $x$ from $\left[-1,1\right]$ distinct from zero, the left-hand side of (\ref{DerivativeEqualZero}) is always positive. Hence, we conclude that the function $f\left(x\right)$ does not have horizontal tangent lines on the interval $\left[-1,1\right]$.

However, a closer look at the graph of the function $f\left(x\right)$ seems to point at a different result: there is a horizontal tangent at $x=0$ (see Figure \ref{fig:FunctionGraph}). 

First, note that the function $f\left(x\right)$ is defined in $x=0$. In order to verify if it has a horizontal tangent at this point, let us find the derivative of the function $f\left(x\right)$ using definition:
\begin{eqnarray}
f'\left(0\right) &=& \lim_{h\rightarrow0}{\frac{f\left(0+h\right)-f\left(0\right)}{h}} \notag \\ 
&=& \lim_{h\rightarrow0}{\frac{\sqrt[3]{h}\sin{\left(h^2\right)}}{h}} \notag \\ 
&=& \lim_{h\rightarrow0}{\left(\sqrt[3]{h} \cdot {h} \cdot \frac{\sin{\left(h^2\right)}}{h^2}\right)} \notag \\
&=& \lim_{h\rightarrow0}{\sqrt[3]{h}} \cdot \lim_{h\rightarrow0}{h} \cdot \lim_{h\rightarrow0}{\frac{\sin{\left(h^2\right)}}{h^2}} \notag \\
&=& 0 \cdot 0 \cdot 1 = 0 \notag
\end{eqnarray}
since each of the limits above exists. We see that, indeed, the function $f\left(x\right)$ possesses a horizontal tangent line at the point $x=0$.

\section{Closer Look at the Expression for the Derivative}

What is the problem with the standard procedure proposed by many textbooks and repeated in every Calculus class? The explanation lies in the following premise: the expression of the derivative of the function does not contain the information as to whether the function is differentiable or not at the points where it is undefined. As it is pointed out in \cite{Rivera2013}, the domain of the derivative is determined \emph{a priori} and therefore should not be obtained from the formula of the derivative itself.

In the example above the Product Law for derivatives requires the existence of the derivatives of both functions at the point of interest. Since the function $\sqrt[3]{x}$ is not differentiable in zero, the Product Rule cannot be applied. 

In order to see what exactly happens when we apply the Product Rule, let us find the expression for the derivative using definition of the derivative:
\begin{eqnarray}
f'\left(x\right) &=& \lim_{h\rightarrow0}{\frac{f\left(x+h\right)-f\left(x\right)}{h}} \notag \\ 
&=& \lim_{h\rightarrow0}{\frac{\sqrt[3]{x+h}\sin{\left(x+h\right)^2}-\sqrt[3]{x}\sin{\left(x^2\right)}}{h}} \notag \\ 
&=& \lim_{h\rightarrow0}{\frac{\left(\sqrt[3]{x+h}-\sqrt[3]{x}\right)}{h}\sin{\left(x^2\right)}} + \notag \\
&& \lim_{h\rightarrow0}{\frac{\left(\sin{\left(x+h\right)^2}-\sin{\left(x^2\right)}\right)}{h}\sqrt[3]{x+h}} \notag \\
&=& \lim_{h\rightarrow0}{\frac{\sqrt[3]{x+h}-\sqrt[3]{x}}{h}} \cdot \lim_{h\rightarrow0}{\sin{\left(x^2\right)}} +  \notag \\&& \lim_{h\rightarrow0}{\frac{\sin{\left(x+h\right)^2}-\sin{\left(x^2\right)}}{h}} \cdot \lim_{h\rightarrow0}{\sqrt[3]{x+h}} \notag \\
&=& \frac{1}{3\sqrt[3]{x^2}} \cdot \sin{\left(x^2\right)}+2x\cos{\left(x^2\right)} \cdot \sqrt[3]{x} \notag 
\end{eqnarray}
which seems to be identical to the expression (\ref{DerivativeExpression}).

Students are expected to develop a skill of deriving similar results and know how to find the derivative of the function using definition of the derivative only. But how `legal' are the performed operations?

\begin{figure}[H]
\begin{center}
	\includegraphics[width=6.0in]{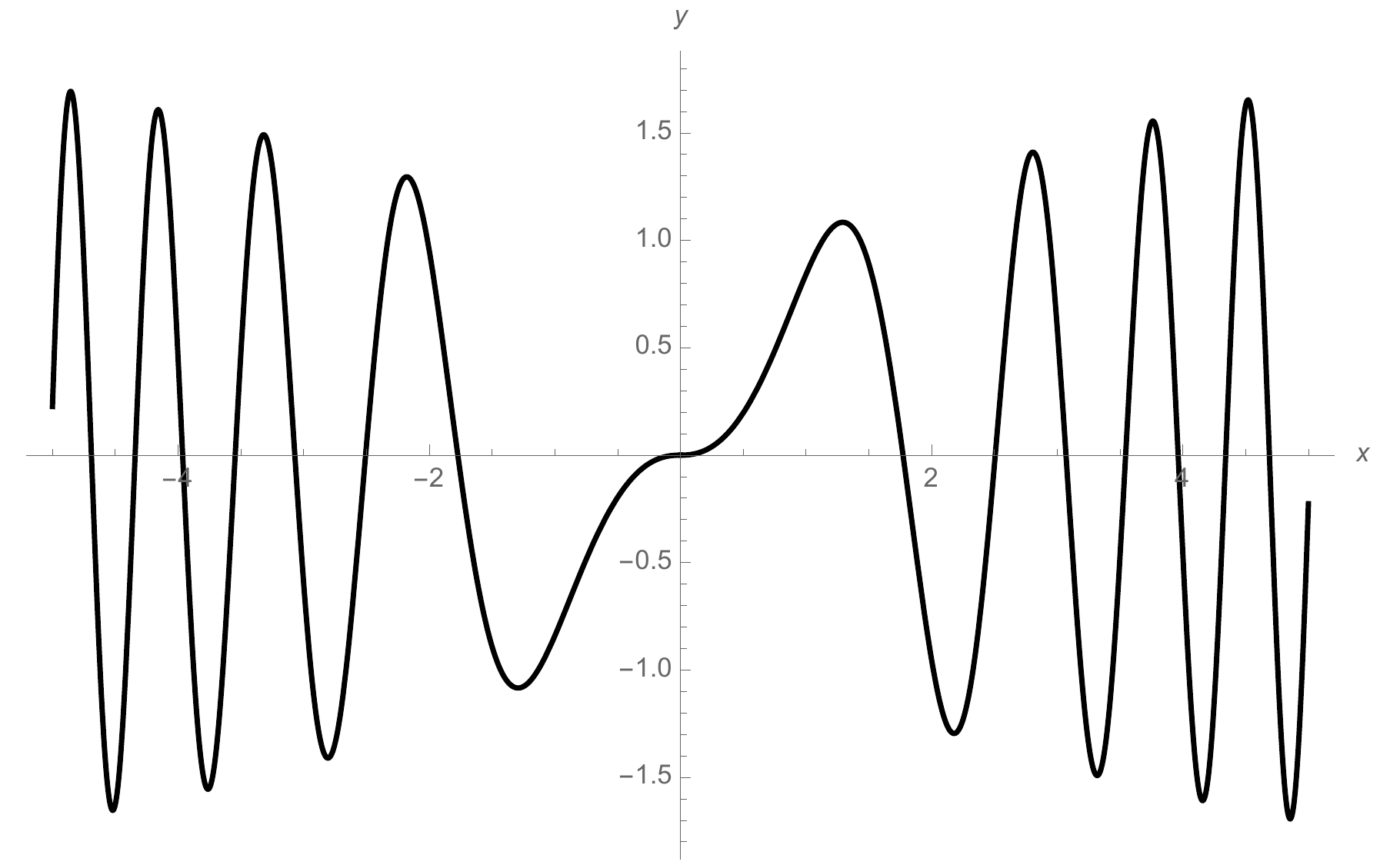}
	\vspace{.1in}
	\caption{Graph of the function $g\left(x\right)=\sqrt[3]{x}\cos{\left(x^2\right)}$}
	\label{fig:GFunction}
\end{center}
\end{figure}

Let us consider each of the following limits: 
\begin{eqnarray*}
&& \lim_{h\rightarrow0}{\frac{\sqrt[3]{x+h}-\sqrt[3]{x}}{h}} \notag \\
&& \lim_{h\rightarrow0}{\sin{\left(x^2\right)}}\notag \\
&& \lim_{h\rightarrow0}{\frac{\sin{\left(x+h\right)^2}-\sin{\left(x^2\right)}}{h}}\notag \\
&& \lim_{h\rightarrow0}{\sqrt[3]{x+h}}.
\end{eqnarray*}
The last three limits exist for all real values of the variable $x$. However, the first limit does not exist when $x=0$. Indeed
\begin{equation*}
\lim_{h\rightarrow0}{\frac{\sqrt[3]{0+h}-\sqrt[3]{0}}{h}} = \lim_{h\rightarrow0}{\frac{1}{\sqrt[3]{h^2}}} = + \infty
\end{equation*}

This implies that the Product and Sum Laws for limits cannot be applied and therefore this step is not justifiable in the case of $x=0$. When the derivation is performed, we automatically assume the conditions, under which the Product Law for limits can be applied, i.e. that both limits that are multiplied exist. It is not hard to see that in our case these conditions are actually equivalent to $x\neq0$. This is precisely why, when we wrote out the expression for the derivative (\ref{DerivativeExpression}), it already contained the assumption that it is only true for the values of $x$ that are different from zero.

Note, that in the case of $x=0$ the application of the Product and Sum Laws for limits is not necessary, since the term $\left(\sqrt[3]{x+h}-\sqrt[3]{x}\right)\sin{\left(x^2\right)}$ vanishes.

The correct expression for the derivative of the function $f\left(x\right)$ should be the following:
\begin{equation*}
f'\left(x\right)  = 
\begin{cases} 
\frac{6x^2\cos{\left(x^2\right)}+\sin{\left(x^2\right)}}{3\sqrt[3]{x^2}}, & \mbox{if } x \neq 0 \\ 
0, & \mbox{if } x = 0 
\end{cases}
\end{equation*}

The expression for the derivative of the function provides the correct value of the derivative only for those values of the independent variable, for which the expression is defined; it does not tell anything about the existence or the value of the derivative, where the expression for the derivative is undefined. Indeed, let us consider the function
\begin{equation*}
g\left(x\right) = {\sqrt[3]{x}}\cos{\left(x^2\right)}
\end{equation*}
and its derivative $g'\left(x\right)$ 
\begin{equation*}
g'\left(x\right) = \frac{\cos{\left(x^2\right)}-6x^2\sin{\left(x^2\right)}}{3\sqrt[3]{x^2}}
\end{equation*}

Similar to the previous example, the expression for the derivative is undefined at $x=0$. Nonetheless, it can be shown that $g\left(x\right)$ is not differentiable at $x=0$ (see Figure \ref{fig:GFunction}). Therefore, we provided two visually similar functions: both have the expressions for their derivatives undefined in zero, however, one of these functions possesses a derivative, but the other one does not.

\section{Methodological Remarks}

Unfortunately, there exist many functions similar to the ones discussed above and they can arise in a variety of typical Calculus problems: finding the points where the tangent line is horizontal, finding an equation of the tangent and normal lines to the curve at the given point, the use of differentials and graph sketching. Relying only on the expression of the derivative for determining its value at the undefined points may lead to missing a solution (as in the example discussed above) or to some completely false interpretations (as in the case of curve sketching).

As it was discussed above, the expression for the derivative does not provide any information on the existence or the value of the derivative, where the expression itself is undefined. Here we present a methodology for the analysis of this type of functions.

Let $f\left(x\right)$ be the function of interest and $f'\left(x\right)$ be the expression of its derivative undefined at some point $x_{0}$. In order to find out if $f\left(x\right)$ is differentiable at $x_{0}$, we suggest to follow a list of steps:

\begin{enumerate}
  \item Check if the function $f\left(x\right)$ itself is defined at the point $x_{0}$. If $f\left(x\right)$ is undefined at $x_{0}$, then it is not differentiable at $x_{0}$. If $f\left(x\right)$ is defined at $x_{0}$, then proceed to next step. 
  \item Identify the basic functions that are used in the formula of the function $f\left(x\right)$, that are themselves defined at the point $x_{0}$, but their derivative is not (such as, for example, the root functions).
	\item Find the derivative of the function $f\left(x\right)$ at the point $x_{0}$ using definition.
\end{enumerate}

The importance of the first step comes from the fact that most students tend to pay little attention to the functions domain analysis when asked to investigate its derivative. Formally, the second step can be skipped, however it will give the students the insight into which part of the function presents a problem and teach them to identify similar cases in the future. the difficulty of accomplishing the third step depends on the form of the function and sometimes can be tedious. Nevertheless, it allows the students to apply the previously obtained skills and encourages the review of the material.

\begin{figure}[H]
\begin{center}
	\includegraphics[width=6.0in]{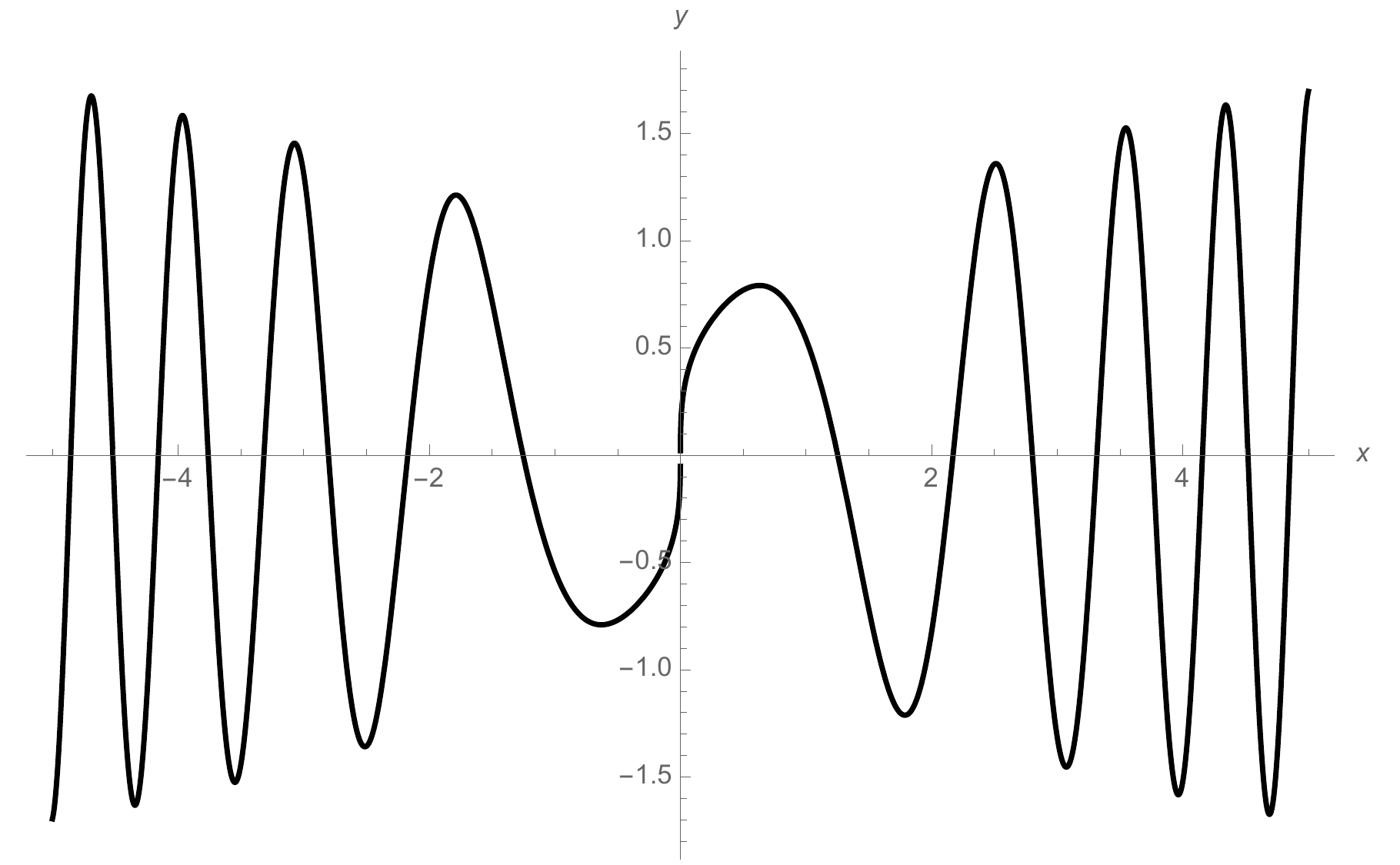}
	\vspace{.1in}
	\caption{Graph of the function $g\left(x\right)=\sqrt[3]{x}\cos{\left(x^2\right)}$}
	\label{fig:GFunction}
\end{center}
\end{figure}

\section{Conclusion}

We discussed the misconception, that the expression of the derivative of the function contains the information as to whether the function is differentiable or not at the points, where the expression is undefined. We considered a typical Calculus problem of looking for the horizontal tangent line of a function as an example. We showed how the search for the values that make the expression of the derivative equal zero leads to missing a solution: even though the expression of the derivative is undefined, the function still possesses the derivative at the point. We provided an example of the function that similarly has the expression for the derivative undefined, however the function is not differentiable at the point. We also presented the methodological treatment of such functions by applying the definition of the derivative, which can be used in the classroom.

\end{document}